



\newif\ifIsbfsINCLUDED
\ifx\IsbfsINCLUDED\undefined\IsbfsINCLUDEDtrue\else\IsbfsINCLUDEDfalse\fi
\ifIsbfsINCLUDED

\documentclass[11pt,leqno,twoside]{article}
\usepackage{amsfonts}
\usepackage{amssymb}
\usepackage{amsthm}

\setlength{\paperheight}{224mm}
\setlength{\paperwidth}{156mm}
\setlength{\textwidth}{120mm}
\setlength{\textheight}{175mm}

\makeatletter
\def\Isbfs@concat#1#2{%
 \expandafter\def\expandafter#1\expandafter{#1#2}}
\def\Isbfs@concat@b#1#2{%
 \expandafter\Isbfs@concat\expandafter#1\expandafter{#2}}

\pagestyle{myheadings}

\def\Isbfs@runningtitle{}
\def\Isbfs@title{}
\def\IsbfsTITLE[#1]#2{\def\Isbfs@runningtitle{#1}\def\Isbfs@title{#2}}

\def\IsbfsAUTHOR[#1]#2{%
\ifx\Isbfs@author@notlast\@empty
\Isbfs@concat@b\Isbfs@author@notlast\Isbfs@author@last
\Isbfs@concat@b\Isbfs@runauthor@notlast\Isbfs@runauthor@last\else
\Isbfs@concat\Isbfs@author@notlast{, }%
\Isbfs@concat@b\Isbfs@author@notlast\Isbfs@author@last
\Isbfs@concat\Isbfs@runauthor@notlast{, }%
\Isbfs@concat@b\Isbfs@runauthor@notlast\Isbfs@runauthor@last\fi
\def\Isbfs@author@last{#2}%
\def\Isbfs@runauthor@last{#1}%
\IsbfsAUTHOR@b#2\\ \@nil}

\def\IsbfsAUTHOR@b#1\\#2\@nil{%
\Isbfs@concat\Isbfs@alladdresses{\medskip\\{\sc #1}\\}}

\def\IsbfsADDRESS#1{\Isbfs@concat\Isbfs@alladdresses{#1\\}}

\def\Isbfs@author@notlast{}
\def\Isbfs@author@last{}
\def\Isbfs@runauthor@notlast{}
\def\Isbfs@runauthor@last{}
\def\Isbfs@alladdresses{}

\def\Isbfs@MAKEauthors{%
 \def\Isbfs@authors{}\def\Isbfs@runauthors{}%
 \ifx\Isbfs@author@notlast\@empty
 \Isbfs@concat@b\Isbfs@author@notlast\Isbfs@author@last
 \Isbfs@concat@b\Isbfs@runauthor@notlast\Isbfs@runauthor@last\else
 \Isbfs@concat\Isbfs@author@notlast{ and }%
 \Isbfs@concat@b\Isbfs@author@notlast\Isbfs@author@last
 \Isbfs@concat\Isbfs@runauthor@notlast{ and }%
 \Isbfs@concat@b\Isbfs@runauthor@notlast\Isbfs@runauthor@last\fi
 \expandafter\Isbfs@MAKEauthors@b\Isbfs@author@notlast\\,\@nil
 \Isbfs@concat@b\Isbfs@runauthors\Isbfs@runauthor@notlast}

\def\Isbfs@MAKEauthors@b#1\\,#2\@nil{%
 \def\Isbfs@authors@b{#2}%
 \ifx\Isbfs@authors@b\@empty
 \Isbfs@concat\Isbfs@authors{#1}\else
 \Isbfs@concat\Isbfs@authors{#1,\\}%
 \expandafter\Isbfs@MAKEauthors@b#2\@nil\fi}

\def\IsbfsEMAIL#1{%
\Isbfs@concat\Isbfs@alladdresses{Email: {\bf #1} \\}}

\def\IsbfsKEYWORDS#1{\def\Isbfs@keywords{#1}}
\def\Isbfs@keywords{}
\def\IsbfsSUBJCLASS#1{\def\Isbfs@subjclass{#1}}
\def\Isbfs@subjclass{}
\def\IsbfsTHANKS#1{\def\Isbfs@thanks{#1}}
\def\Isbfs@thanks{}
\def\IsbfsDEDICATEDTO#1{\def\Isbfs@dedicatedto{#1}}
\def\Isbfs@dedicatedto{}

\long\def\IsbfsABSTRACT#1{\long\def\Isbfs@abstract{#1}}

\def\IsbfsMAKETITLE{
\thispagestyle{plain}
\Isbfs@MAKEauthors
\markboth{\Isbfs@runauthors}{\Isbfs@runningtitle}

\begin{center}
\large\bf
\Isbfs@title
\end{center}

\begin{center}
\Isbfs@authors
\end{center}
\Isbfs@definemacros

\ifx\Isbfs@dedicatedto\@empty\smallskip\else
\begin{center}\Isbfs@dedicatedto\end{center}\fi

\begin{quote}
{\sc Abstract.}
\Isbfs@abstract
\ifx\Isbfs@subjclass\@empty\relax\else
{\def\thefootnote{{}}%
\footnote{2000 {\it subject classification}: \Isbfs@subjclass}%
\addtocounter{footnote}{-1}}\fi
\ifx\Isbfs@keywords\@empty\relax\else
{\def\thefootnote{{}}%
\footnote{{\it Key words and phrases}: \Isbfs@keywords}%
\addtocounter{footnote}{-1}}\fi
\ifx\Isbfs@thanks\@empty\relax\else
{\def\thefootnote{{}}%
\footnote{\Isbfs@thanks}%
\addtocounter{footnote}{-1}}\fi
\end{quote}
}

\long\def\IsbfsDEFINEMACROS#1{\long\def\Isbfs@definemacros{#1}}

\def\IsbfsBIBLIOGRAPHY{
\begin{thebibliography}{99}
\itemsep=6pt plus 1pt}

\def\endIsbfsBIBLIOGRAPHY{
\end{thebibliography}}

\def\Isbfs@newline{\par\noindent}
\def\endIsbfsDOCUMENT{%
\begingroup\let\\=\Isbfs@newline\Isbfs@alladdresses\endgroup
\enddocument}

\newcounter{Suzuki@citeA}
\newcounter{Suzuki@citeB}
\newcounter{Suzuki@citeC}
\newcounter{Suzuki@citeD}
\def\IsbfsCITE#1{%
 \begingroup[\setcounter{Suzuki@citeA}{-1}\Suzuki@Cite{#1}]\endgroup}
\def\Suzuki@Cite#1{%
 \let\Suzuki@Cite@a\@empty
 \let\Suzuki@Cite@b\@empty
 \@for\Suzuki@Cite@i:=#1\do{%
 \edef\Suzuki@Cite@i{\expandafter\@firstofone\Suzuki@Cite@i\@empty}%
 \@ifundefined{b@\Suzuki@Cite@i}{%
  \Suzuki@Cite@b\let\Suzuki@Cite@b\@empty
  \Suzuki@Cite@a
  \setcounter{Suzuki@citeA}{-1}%
  \mbox{\reset@font\bfseries ?}%
  \G@refundefinedtrue
  \@latex@warning{Citation `\Suzuki@Cite@i' on page \thepage \space
    undefined}}%
 {\setcounter{Suzuki@citeB}{\csname b@\Suzuki@Cite@i\endcsname}%
  \ifnum\value{Suzuki@citeA}=\value{Suzuki@citeB}%
   \setcounter{Suzuki@citeD}{\the\value{Suzuki@citeB}}%
   \ifnum\value{Suzuki@citeC}<\value{Suzuki@citeB}%
    \def\Suzuki@Cite@b{--\Suzuki@Cite@Form{\the\value{Suzuki@citeD}}}%
   \else\def\Suzuki@Cite@b{,\penalty\@m
    \ \Suzuki@Cite@Form{\the\value{Suzuki@citeD}}}%
   \fi\addtocounter{Suzuki@citeA}{1}%
  \else
   \Suzuki@Cite@b\let\Suzuki@Cite@b\@empty
   \Suzuki@Cite@a
   \hbox{\Suzuki@Cite@Form{\the\value{Suzuki@citeB}}}%
   \setcounter{Suzuki@citeA}{\the\value{Suzuki@citeB}}%
   \addtocounter{Suzuki@citeA}{1}%
   \setcounter{Suzuki@citeC}{\the\value{Suzuki@citeA}}%
  \fi}%
 \def\Suzuki@Cite@a{,\penalty\@m\ }}%
 \Suzuki@Cite@b}
\def\Suzuki@Cite@Form#1{#1}

\def\IsbfsINDENTAFTERSECTION{
\IsbfsINDENTAFTERSECTION@b\section\Isbfs@SAVE@SECTION
\IsbfsINDENTAFTERSECTION@b\subsection\Isbfs@SAVE@SUBSECTION
\IsbfsINDENTAFTERSECTION@b\subsubsection\Isbfs@SAVE@SUBSUBSECTION
}

\def\IsbfsINDENTAFTERSECTION@b#1#2{
\let#2#1
\def#1{\@ifstar{\IsbfsINDENTAFTERSECTION@c#2}%
{\IsbfsINDENTAFTERSECTION@d#2}}}

\def\IsbfsINDENTAFTERSECTION@c#1#2{#1*{#2}\@afterindenttrue}

\def\IsbfsINDENTAFTERSECTION@d#1#2{#1{#2}\@afterindenttrue}

\makeatother
\fi



\IsbfsTITLE[running title]{Compactness and Measures of
Noncompactness in \\ Metric Trees }

\IsbfsAUTHOR[A. G. Aksoy]{A. G. Aksoy} \IsbfsADDRESS{Department of
Mathematics, Claremont McKenna College, Claremont, CA  91711, USA. }
\IsbfsEMAIL{aaksoy@cmc.edu}

\IsbfsAUTHOR[M. S. Borman]{M. S. Borman}
\IsbfsADDRESS{Department of Mathematics, University of Chicago, Chicago, IL
60637, USA} \IsbfsEMAIL{borman@math.uchicago.edu}

\IsbfsAUTHOR[A. Westfahl]{A. L. Westfahl} \IsbfsADDRESS{Department
of Mathematics, Claremont McKenna College, Claremont, CA  91711,
USA.} \IsbfsEMAIL{awestfahl07@cmc.edu}

\IsbfsTHANKS{Research supported by  the NSF-REU grant DMS 0453284}

\IsbfsDEFINEMACROS{

\newtheorem{theorem}{Theorem}[section]
\newtheorem{lemma}[theorem]{Lemma}
\newtheorem{corollary}[theorem]{Corollary}

\theoremstyle{definition}
\newtheorem{definition}[theorem]{Definition}
\newtheorem{example}[theorem]{Example}

\theoremstyle{remark}
\newtheorem{remark}[theorem]{Remark}



}


\newcommand{\s}[2]{S_{#1,#2}}
\newcommand{\R}{\mathbb{R}}

\newcommand{\N}{\mathbb{N}}
\newcommand{\iso}[2]{h_{#1,#2}}

\newcommand{\ep}{\epsilon}



\IsbfsSUBJCLASS{Primary 54E45, Secondary 05C05, 47H09}
\IsbfsKEYWORDS{Metric Trees, Measures of Noncompactness, Condensing
Operators}

\IsbfsABSTRACT{

A metric tree ($M$, $d$), also known as $\mathbb{R}$-trees or
$T$-theory, is a metric space such that between any two points there
is an unique arc and that arc is isometric to an interval in
$\mathbb{R}$.

In this paper after presenting some fundamental properties of metric
trees and metric segments, we will give a characterization of
compact metric trees in terms of metric segments. Two common
measures of noncompactness are the ball and set measures,
respectively defined as
\[
\beta(S) = \inf \{ \epsilon \mid S \mbox { has a finite }
\epsilon\mbox {-net in } M\} \mbox{ and }
\]
\[
\alpha(S) = \inf\{\epsilon \mid S \mbox{ has a finite cover of sets
of diameter} \leq \epsilon\}.
\]
 We will prove that  $\alpha = 2\beta$ for all metric trees.
 We give two independent proofs of this result, first of which depends
 on the fact that given a metric tree $M$ and a subset $E$ with its diameter $2r$,
 for all $\ep >0$ there exists $m \in M$ such that $E\subset~\hspace*{-3pt}B(m;r+\ep)$.
 The second proof depends on the  calculation of underlying geometric constant,
   Lifschitz characteristic of a metric tree.
 It is well known that the classes of condensing or contractive
 operators defined relative to distinct measures of noncompactness
 are not equal in general, however in case of metric trees, we show
 that a map $T$ between two metric trees is $k$-set contractive if
 and only if it is $k$-ball contractive for $k\geq 0$.
}

\begin{IsbfsDOCUMENT}

\section{Introduction}

The study of metric trees (T-theory or $\R$-trees) began with J.
Tits \cite{Tits} in 1977 and since then, applications have been
found for metric trees within many fields of mathematics. For an
overview of geometry, topology, and group theory applications,
consult Bestvina \cite{Best}. A complete discussion of these spaces
and their relation to so-called $CAT(0)$ spaces we refer to
\cite{Brid}. Applications of metric trees in biology and medicine
involve phylogenetic trees \cite{Semp}, and metric trees even find
applications in computer science that involve string matching
\cite{Bart}.  Given that metric trees go by three names there are a
couple of equivalent definitions, and to understand the definition
we choose to use we must first define the notion of a metric segment
in a metric space. 
\begin{definition}\label{D:mseg}
    Let $x, y \in X$, where ($X$, $d$) is a metric space.
    A \emph{metric segment from $x$ to $y$}, denoted
    $[x,y]$, is a subset of $X$ such that $\iso{x}{y}:z \mapsto d(x,z)$ is an
    isometry from $[x, y]$ onto $[0, d(x,y)] \subset \R$ where $\R$
    has the Euclidian metric. 
\end{definition}

Following the notation used for segments in $\R$, for
a metric segment $[x,y] \subset X$ denote a half-open segment
as $[x,y) := [x, y] \setminus \{y\}$ and denote an open segment from 
$x$ to $y$ as $(x, y) := [x,y] \setminus \{x, y\}$.

\begin{definition}\label{D:mt2}
    $(M,d)$, a metric space, is a \emph{metric tree} iff for all $x,y,z\in M$,
    \begin{enumerate}
        \item there exists an unique metric segment from $x$ to $y$,
\vspace*{-4pt}
        \item $[x,z] \cap [z,y] = \{z\} \Rightarrow [x,z] \cup [z,y] = [x,y]$.
    \end{enumerate}
\end{definition}
Throughout this paper, ($X$, $d$) will be a
metric space and ($M$, $d$) will be a metric tree.  Furthermore,
$d(x,y)$ will be denoted as $xy$ as long as there is no fear of
confusion, i.e., \vspace*{-4pt}
$$ d(x,y):= xy.
\vspace*{-4pt}$$ 
Let $B(x;r):=\{y \in X \mid xy < r\}$ denote the open ball centered at $x$ with
radius $r$ and $B_c(x;r) := \{y \in X \mid xy \leq r\}$ denote the closed
ball centered at $x$ with radius $r$.

\vspace*{-4pt}
\begin{example}[A Nice Simple Tree]
    Let $M$ be the subset of $\R^2$ where
    $M = [-1,1] \times \{0\} \; \bigcup \; \{1\} \times [-1,1]$.
    Define a metric $d: M \times M \to \R$ by:

           \[ d(x,y) = \left\{\begin{array}{ll}
                              \|x-y\| & \mbox{if  $x,y \in [-1,1] \times \{0\}$,}\\
                              \|x-y\| & \mbox {if $x,y \in \{1\} \times [-1,1]$,}\\
                              \|x-(1,0)\| + \|(1,0)-y\| & \mbox{otherwise.}
                              \end{array}
                       \right.\]

    With this metric, ($M$, $d$) becomes a metric tree.
    While $M$ looks like a tree from graph theory, this is not the case
    for all metric trees.
\end{example}

\begin{example}[The Radial Metric]\label{E:radial}
    Define $d: \R^2 \times \R^2 \to \R$ by:

        \[ d(x,y) =\left \{\begin{array}{ll}
                        \|x-y\| & \mbox {if $x = \lambda \, y$ for some $\lambda \in \R$,}\\
                        \|x\| + \| y \| & \mbox{otherwise.}
                        \end{array}
                 \right.\]
    We can observe the $d$ is in fact a metric and that $(\R^2,d)$ is a metric tree.

    In order to investigate the radial metric, we will actually investigate a metric tree
    that is isometric to the radial metric.  Let $M$ be defined as:
    \[
    M:= \left\{x \in \prod_{\alpha \in [0, 2\pi)}\R_{\geq 0} \mid
    x_\alpha \not= 0 \mbox{ for at most one } \alpha \right\}.
    \]
    Then, define a metric on $M$ as:
    \[
        \rho(x,y) = \sum_\alpha |x_\alpha - y_\alpha|.
    \]
    The summation is defined, for $|x_\alpha - y_\alpha| \not= 0$ for at
    most two $\alpha \in [0, 2\pi)$.
    $(\R^2,d)$ and $(M,\rho)$ are isometric by the isometry $f: \R^2 \to M$ with
    \[ f(r\,e^{i\beta})_\alpha = \left\{\begin{array}{ll}
                                          r & \mbox{if $\alpha = \beta$,}\\
                                          0 & \mbox{otherwise.}
                                          \end{array}
                                  \right.\]

    Let $\mathbf{0}$ denote the element of $M$ that is zero in every coordinate.
    Observe that $f^{-1}(B(\mathbf{0};r)_\rho) = B((0,0);r)_{Euclid}$.
    For $x \in M$ with $x_\beta \not= 0$, (here $f(x_\beta \, e^{i\beta}) = x$),
    \[
        B(x;r)_\rho = f(\{t\,e^{i\beta}\mid t \in (x_\beta - r, x_\beta + r)\})
        \cup B(\mathbf{0};r-x_\beta)_\rho
    \]

    Suppose $x^n \rightarrow x$ in $(M,\rho)$ with $x_\beta \not= 0$.
    Then, there exists $N \in \N$ such that
    $n \geq N \Rightarrow x^n \in B(x;x_\beta)_\rho$.
    Hence $n \geq N \Rightarrow x^{n}_{\alpha}=0$ for $\alpha \not= \beta$ and
    $x^{n}_\beta \in (0,2x_\beta)$,
    since $B(x;x_\beta)_\rho = f(\{t\,e^{i\beta}\mid t \in (0, 2x_\beta)\})$.
    In terms of the radial metric, this means that the tail end of $\{x^n\}$ is a
    convergent sequence in $[0,2x_\beta] \subset \R$.
\end{example}

\subsection{Properties of Metric Segments}
This subsection covers some basic properties of metric segments in
metric spaces drawn from Blumenthal \cite{Blum}.  Let ($X$, $d$) be
a metric space and for $a,b \in X$ let $\s{a}{b}$ denote a metric
segment from $a$ to $b$ (since metric segments are not necessarily
unique).

\begin{definition}\label{D:betw}
    Let $x, y, z \in X$.  $y$ is \emph{between $x$ and $z$}, denoted
    $xyz$, iff $xz = xy + yz$.
\end{definition}

\begin{remark}
    Given metric segment $\s{x}{z} \subset X$, observe that $xyz$ for all $y \in \s{x}{z}$.
\end{remark}

\begin{lemma}[Blumenthal, \cite{Blum}]\label{L:trans}
    Let $a,b,c,d \in X.$ If $abc$ and $acd$, then $abd$ and $bcd$.
\end{lemma}

 We will often refer to above property of metric segments, which
is about the transitivity of betweenness of points on metric
segments.
\begin{lemma}[Blumenthal, \cite{Blum}]\label{L:B151}
    Let $a,b,p \in X$ and suppose $\s{a}{p}$ and $\s{p}{b}$ exist in $X$.
    $\s{a}{p} \cup \s{p}{b} = \s{a}{b}$ iff $apb$.
\end{lemma}
\begin{theorem}[Blumenthal, \cite{Blum}]\label{T:B151}
    Let $X$ be a complete metric space, $F$ a closed subset of $X$, $a,b \in F$,
    and $F \subset \{x \in X \mid axb \}$.
    Then $F$ is a metric segment from $a$ to $b$ iff
    $F$ has exactly one midpoint of $x,y$ for all $x,y \in F$.
\end{theorem}
\begin{theorem}[Blumenthal, \cite{Blum}]\label{T:B152}
    Let $a,b \in X$ and let $F \subset X$ be an injective arc from
    $a$ to $b$.  $F$ is a metric segment iff for all $x,y \in F$ either $axy$ or $yxb$.
\end{theorem}


\subsection{Metric Segments in Metric Trees}

Within this section ($M$, $d$) will be a metric tree.

\begin{lemma}\label{L:ms} ~
\begin{enumerate}
\item  For $x,y \in M$, $[x,y] = \{z \in M \mid xzy\}$.
\item For $x,y,z \in M$, $z \in [x,y]$ only if $[x,z] \subset [x,y]$.
\end{enumerate}
\end{lemma}
The proof is classical and follows directly from the results in
\cite{Blum}.

\begin{theorem}\label{T:msip}
    If $x,y,z \in M$, then there exists $w \in M$ such that $[x,z] \cap [y,z] = [w,z]$
    and $[x,y] \cap [w,z] = \{w\}$.
\end{theorem}
\begin{proof}
    Define $f:[x,z] \cap [y,z] \to \R$ by $v \mapsto vz$.
    Since $[x,z] \cap [y,z]$ is compact, there exists $w \in [x,z] \cap [y,z]$
    such that $wz \geq vz$ for all $v \in [x,z] \cap [y,z]$.

    By Lemma ~\ref{L:ms}, $[w,z] \subset [x,z] \cap [y,z]$.
    Conversely let $a \in [x,z] \cap [y,z]$, then we know
    that either $xaw$ or $waz$ by Theorem~\ref{T:B152}.
    If $xaw$, then by transitivity and $xwz$ we know that $awz$.
    Hence $aw + wz = az \leq wz$ using the maximality of $wz$,
    hence $aw = 0$ and $a = w$.  Hence $a = w$ or $waz$, and in either case
    $a \in [w,z]$ by Lemma~\ref{L:ms}.  Therefore $[x,z] \cap [y,z]=~[w,z]$.

    Suppose $[x,w] \cap [w,y] \not= \{w\}$.
    Let $w' \in [x,w] \cap [w,y] \subset [x,z] \cap [y,z]$ with $w' \not= w$.
    Hence $xw'w$ and $xwz$, and therefore by transitivity $xw'z$.
    By $xw'z$, the triangle inequality, $xw'w$, and $xwz$,
    \[
    xz = xw' + w'z \leq xw' + w'w + wz = xw + wz = xz.
    \]
    Therefore, $w'z = w'w +wz > wz$, which contradicts the maximality of $wz$.

    Hence $[x,w] \cap [w,y] = \{w\}$, therefore $[x,w] \cup [w,y] = [x,y]$ and
    in particular $w \in [x,y]$.
    Let $a \in [x,y] \cap [w,z]$ and $a \not= w$. Since $a, w \in [x,y]$
    by Theorem~\ref{T:B152}, we know without loss
    of generality that $xaw$.  By transitivity, $xaw$ and $xwz$ implies $awz$.
    This means that $az = aw + wz > wz$, which contradicts the maximality of $wz$.
    Therefore, $[x,y] \cap [w,z]=~\{w\}$.
\end{proof}

\begin{theorem}\label{T:midpointball}
    Let $a,m,x,y \in M$.
    Suppose $m$ is the midpoint of $[x,y]$ and $x,y \in B(a;r)$.
    Then $B_c(m;\frac{xy}{2}) \subset B(a;r)$.
\end{theorem}
\begin{proof}
    Let $z \in B_c(m;\frac{xy}{2})$, so $zm \leq \frac{xy}{2}$,
    we want to show that $za < r$.

    By Theorem~\ref{T:msip}, exists $w \in [x,y]$ such that $[a,x] \cap [a,y] = [a,w]$.
    Since $w \in [x,y]$, by Theorem~\ref{T:B152} we can assume
    without loss of generality that $xmw$.  By transitivity, $xwa$ and $xmw$ imply
    $xma$.  So $\frac{xy}{2} + ma = xm + ma = xa < r$ and hence $ma< r-\frac{xy}{2}$.
    By the above, $za \leq zm + ma < \frac{xy}{2} + (r-\frac{xy}{2}) = r$.

    Hence $za < r$ and therefore $B_c(m;\frac{xy}{2}) \subset B(a;r)$.
\end{proof}

\begin{corollary}\label{C:bc}
    If $x,y,a \in M$ with $x,y \in B(a;r)$, then $[x,y] \subset B(a;r)$.
\end{corollary}
\begin{proof}
    Let $m$ be the midpoint of $[x,y]$.
    Then $[x,y] \subset B_c(m;\frac{xy}{2})$ and
    therefore by Theorem~\ref{T:midpointball}, $[x,y] \subset B(a;r)$.
\end{proof}
\section{Compactness in Metric Trees}

Since metric trees are very similar to $\R$, one might hope that the
Heine-Borel theorem is true in metric trees.  Unfortunately this is
not the case as the following example demonstrates.

\begin{example}\label{E:unitradial}
    Returning to Example~\ref{E:radial} of the radial metric on $\R^2$,
    consider the closed unit ball at the origin $B_c((0,0);1)$.  For 
    $\alpha, \beta \in [0, 2\pi)$, we know that $e^{i\alpha} \in B_c((0,0);1)$ and
    $d(e^{i\alpha},e^{i\beta}) = 2$ for $\alpha \not= \beta$, which implies
    that $B_c((0,0);1)$ is not totally bounded.  Therefore, $B_c((0,0);1)$
    is not compact.
\end{example}

\begin{example} In the following we consider the river metric and
define a compact and noncompact tree with respect to this metric.
These simple examples motivate the characterization of compactness
in metric trees.
\noindent
Let $M=[0,1]^2$ and define a river metric
$\rho:M\times M\rightarrow \R_{\geq 0}$ by:

\[ \rho((x_1,y_1),(x_2,y_2))=\left \{ \begin{array}{ll}
                                    \vert y_1-y_2 \vert & \mbox{$x_1=x_2$}\\
                                     y_1+\vert x_1-x_2\vert +y_2 & \mbox{otherwise}
                                     \end{array}
                              \right. \]
\noindent
Then \emph{a compact tree} $C$ looks like:
 $\displaystyle C=[0,1] \times \{0\} \cup
\bigcup^{\infty}_{n=1} \{\frac{1}{n}\}\times [0,\frac{1}{n}]$

\noindent
and \emph{a noncompact tree} $A$ is: $\displaystyle A=[0,1] \times
\{0\} \cup \bigcup ^{ \infty}_{n=1} \{\frac{1}{n}\} \times [0,1]$.
\end{example}
\begin{definition}\label{D:final}
    Let $M$ be a metric tree.  Define $F$, the set of \emph{final points of $M$(or leaves of $M$)} as
    \[
    F:= \{f \in M \mid f \notin (x,y) \mbox{ for all } x,y \in M\}.
    \]
\end{definition}

\begin{theorem}\label{T:final cover}
    If $M$ be a compact metric tree and $a \in M$, then
    \[
    M = \bigcup_{f \in F} [a,f]
    \]
\end{theorem}
\begin{proof}
    It suffices to show that $M \subset \bigcup_{f \in F} [a,f]$.
    Let $m \in M$, define
    \[
    R_m = \{z \in M \mid m \in [a,z]\}
    \]
    and observe that $R_m \not= \emptyset$ since $m \in R_m$.

    We will show that $R_m$ is closed in $M$.  Let $y \notin R_m$.
    Hence there exists $u \in [a,m)$ such that $[a,y] \cap [m,y] = [u,y]$.
    Let $\ep < um$ and let $b \in B(y; \ep)$.  Hence,
    $$
      \begin{array}{rcll}
      ab&\leq& ay+yb < ay + um&\mbox{since $yb < \ep < um$}\\
        &=& au + uy + um = am + uy&\mbox{since $auy$ and $aum$}\\
        &=& am + (my - mu)&\mbox{since $muy$}\\
        &<& am + my - yb&\mbox{since $yb < \ep < um$}\\
        &\leq& am + mb&\mbox{since $my \leq mb + by$},
    \end{array}  
    $$    
    and therefore $ab < am + mb$, so $b \notin R_m$.
    Hence $B(y;\ep) \subset M \setminus R_m$, and therefore $R_m$ is
    closed in $M$.  Since $M$ is compact it follows that $R_m$ is compact.

    Define $h: R_m \to \R$ by $z \mapsto az$.  Since $R_m$ is compact, there exists
    $f \in R_m$ such that $af \geq az$ for all $z \in R_m$.

    We will show that $f \in F$.  Suppose that $f \notin F$, so there exists $x,y \in M$
    such that $f \in (x,y)$.  Since metric segments are closed under intersections,
    $f \in (x,y)$, and $f$ is an end point of $[a, f]$, we have that
    $[a, f] \cap [x, y] = [f, v]$ where $v \in [x, y]$.
    Hence by switching $x$ and $y$ we can claim without loss of generality that
    $v \in [a, f]$.  Therefore, $f \in [a,y]$ and since $m \in [a,f]$, by transitivity we
    have that $m \in [a,y]$.  Hence $y \in R_m$ and $ay = af + fy > af$, since
    $f \not= y$.  We have a contradiction with the maximally of $af$ in
    $R_m$, therefore $f \in F$.

    Hence for all $m \in M$ there exists $f \in F$ such that $m \in [a,f]$ and
    therefore we have proven that $M \subset \bigcup_{f \in F} [a,f]$.
\end{proof}
\begin{theorem}\label{T:final compact}
    A metric tree $M$ is compact if and only if
    \begin{enumerate}
        \item $M = \bigcup_{f \in F} [a,f]$ for all $a \in M$, and
        \item $\overline{F}$ is compact.
    \end{enumerate}
\end{theorem}
\begin{proof}
    That $M$ is compact implies (1) and (2) follows from Theorem~\ref{T:final cover}, so we just
    need to prove that (1) and (2) imply that $M$ is compact.

    Let $\mathfrak{U}=\{ B(x_\alpha; r_\alpha)\}_\alpha$ be an open-ball cover of $M$.
    Since $\overline{F}$ is compact, we know that there exists a finite subcover of $\overline{F}$
    in $\mathfrak{U}$.  Let $\bigcup_{j=1}^{n} B(x_j, r_j) \supset \overline{F}$ be such a subcover.
    If $M \subset  \bigcup_{j=1}^{n} B(x_j, r_j)$, then we are done, so suppose there exists
    $a \in M$ such that $a \notin \bigcup_{j=1}^{n} B(x_j, r_j)$.

    Let $f, f' \in F \cap B(x_j,r_j)$. Recall that there exists $c \in [f, f']$ such that
     $[a , f] \cap [a, f'] = [a, c]$ and that $[f,f'] \subset B(x_j,r_j)$.  It follows then that
    $[a, f] \setminus B(x_j,r_j) = [a, f'] \setminus B(x_j,r_j)$.
    
    For each $j$, let $f_j \in F \cap B(x_j,r_j)$. It follows from the above that
    $$\displaystyle \bigcup_{f\in F\cap B(x_j ,r_j)} [a , f]=[a ,f_j]\cup B(x_j, r_j).$$  Then by (1) and the above claim we have that
    \[
    M = \bigcup_{f \in F} [a,f] = \bigcup_{j=1}^n [a, f_j]\cup B(x_j,r_j) = \left(\bigcup_{j=1}^n [a, f_j]\right)
    \cup \left(\bigcup_{j=1}^n B(x_j,r_j)\right).
    \]
    Since each $[a, f_j]$ is compact, we can find a finite subcover for it in $\mathfrak{U}$, call it
    $\mathfrak{U}_j$.  Then we have that $\{B(x_j,r_j) \mid 1\leq j \leq n\} \cup
    \bigcup_{j=1}^n \mathfrak{U}_j$ is a finite subcover of $M$ in $\mathfrak{U}$.  Therefore, $M$
    is compact and the theorem is proven.
\end{proof}
\section{Measures of Noncompactness in Metric Trees}

The notion of the measure of noncompactness of a subset of a metric
space was introduced by Kuratowski \cite{Kura} as a way to
generalize Cantor's intersection theorem.  In 1955, Darbo
\cite{Darb} applied measures of noncompactness to prove a powerful
fixed point theorem and since then measures of noncompactness have
been a standard notion in fixed point theory.
A\emph{set-contraction} is a mapping under which the image of any
set is, in some definite sense, more compact than the set itself.
The class of set (or ball)-contractions (or condensing) mappings has
many applications in the study of nonlinear operators. This class is
defined via the notion of measure of noncompactness of a set (mnc
for brevity). A continuous map $T$ is said to be \emph{k-set
contraction} if there is a number $k\geq 0$ such that $$ \alpha
(T(D))\leq k \alpha(D)$$ for all bounded sets $D$. Every compact $T$
is a k-set contraction with $k=0$. Compactness plays an essential
role in the proof of Schauders fixed point theorem, however there
are some situations where the operators are not compact. Darbo
\cite{Darb} utilized the concept of set-contractions by proving a
generalized version of the Schauders fixed point theorem. Besides
its use in fixed point theory, mncs are also a useful tool in
studying theory of functional equations, including partial
differential and integral equations, optimal control theory, etc.
Although all measures of noncompactness are equivalent in a
topological sense, the contractive constants of the mappings are not
preserved when different measures are considered. Moreover, these
measures are closely related to the geometrical properties of the
underlying space. Therefore it is important to know the
relationships between different mncs. In the following we introduce
several mncs in the context of metric trees and study the
relationship between them. For fixed point theorems for metric trees
we refer to \cite{Akso}, \cite{Kirk1} and \cite{Kirk2}.

\begin{definition}\label{D:smnc}
    Given metric space $X$ and bounded subspace $A$, define the
    \emph{set (Kuratowski) measure of noncompactness} as
    \[
    \alpha (A) := \inf\,\left\{b>0 \mid A \subset \bigcup_{j=1}^{n} E_j
    \mbox{ for some } E_j \subset A,\, diam(E_j) \leq b \right\}.
    \]
\end{definition}
\begin{definition}\label{D:bmnc}
    Given metric space $X$ and bounded subset $A$, define the
    \emph{ball (Hausdorff) measure of noncompactness} as
    \[
    \beta (A) := \inf\,\left\{b>0 \mid A \subset \bigcup_{j=1}^{n} B(x_j;b)
    \mbox{ for some } x_j \in X \right\}.
    \]
\end{definition}

\begin{lemma}[Banas and Goebel, \cite{Bana}]\label{L:ba2b}
    For a bounded set $A$ of a metric space $X$,
    \[
    \beta(A) \leq \alpha(A) \leq 2\beta(A).
    \]
\end{lemma}
Before proceeding to prove statements about measures of
noncompactness we will   point out the following  behavior of
$\alpha$ and $\beta$ under isometries. We omit the straightforward
proofs.

\begin{theorem}\label{T:iaie}
    If $(X,d)$ and $(Y,\xi)$ are metric spaces, $f: X \rightarrow Y$ an
    isometric embedding, and $A \subset X$ bounded, then $\alpha(f(A))=\alpha(A)$.
\end{theorem}

\begin{lemma}
    If $(X,d)$ and $(Y,\xi)$ are metric spaces,
    $f: X \to Y$ is an isometric embedding, and $A \subset X$ is bounded,
    then it is not necessarily the case that $\beta(f(A)) = \beta(A)$.
\end{lemma}
\begin{proof}
    Let ($X$, $d$) be a discrete metric space with an infinite number of points.
    Since $X = B(x;1+\ep)$ for all $\ep >0$ and $\{x\} = B(x;1)$, we know
    $\beta(X) = 1$.

    $Y = X \cup \{x^\star\}$ where $x^\star \notin X$.  Make ($Y$, $\xi$) a metric
    space by defining $\xi$ by
    \[
        \xi(a,b) =\left \{\begin{array}{ll}
                            1 & \mbox{if $a,b \in X$}\\
                            \frac{1}{2} & \mbox{if exactly one of $a,b$ is $x^{\star}$}
        \end{array}
        \right.
    \]

    Consider the isometric embedding $f: X \to Y$ by $f(x)=x$.
    Since $f(X) \subset B(x^\star;\frac{3}{4})$, we know that
    $\beta(f(X)) < 1 = \beta(X)$.
    Therefore, $\beta(f(X)) \not= \beta(X)$.
\end{proof}

By the above Lemma we see that $\beta$ is not invariant under
isometric embeddings, however, it is not difficult to prove $\beta$
is invariant under bijective isometries.

We will now turn our attention back to measures of noncompactness in
metric trees.
\begin{lemma}\label{L:bcinmt}
    Let $E$ be a subset of metric tree $M$ with $diam(E) = 2r$.
    Then for all $\ep >0$ there exists $m \in M$ such that $E \subset B(m;r+\ep)$.
\end{lemma}
\begin{proof}
    For all $\ep >0$, there exists $x,y \in E$ such that $xy > 2r - 2\ep$ and let
    $m$ be the midpoint of $[x,y]$.

    Let $z \in E$, then by Theorem~\ref{T:msip},
    there exists $w \in [x,y]$ such that $[z,x] \cap [z,y] = [z,w]$.  Without loss
    of generality we can assume that $m \in [w,x]$ and hence $w \in [z,x]$
    by transitivity.  Next, $diam(E) = 2r$ and $w \in [z,x]$ imply that
    \[
    2r \geq zx = zm + mx = zm + \frac{xy}{2} > zm + (r-\ep),
    \]
    which implies $r + \ep > zm$.  Therefore $E \subset B(m; r+\ep)$.
\end{proof}

\begin{theorem}\label{T:a=2bmt}
    Let $A$ be a bounded subset of metric tree $M$.  Then $\alpha(A) = 2\beta(A)$.
\end{theorem}
\begin{proof}
    By Lemma~\ref{L:ba2b} it suffices to prove that $\alpha(A) \geq 2\beta(A)$.

    Let $2r > \alpha(A)$, then there exists $E_j \subset M$ such that
    $diam(E_j) \leq 2r$ and $A \subset \bigcup_{j=1}^{n} E_j$.
    By Lemma~\ref{L:bcinmt}, for all $\ep >0$ there exists $m_j \in E_j$ such that
    $E_j \subset B(m_j; r+\ep)$.  Hence $A \subset \bigcup_{j=1}^{n} B(m_j;r+\ep)$,
    so $\beta(A) \leq r +\ep$ for all $\ep >0$.  Therefore, $\beta(A) \leq r$.

    So $\beta(A) \leq r$ for all $r$ such that $2r > \alpha(A)$.
    Therefore $\alpha(A) \geq 2\beta(A)$.
\end{proof}

\begin{remark}\label{R:beta*}
    Some authors define the \emph{ball measure of noncompactness} as:
    \[
    \beta^* (A) := \inf\,\left\{b>0 \mid A \subset \bigcup_{j=1}^{n} B(x_j;r_j)
    \mbox{ with }
    diam(B(x_j,r_j)) \leq b \right\}.
    \]
\end{remark}

\begin{theorem}\label{T:b* = b mt}
    If $A$ be a bounded subset of a metric tree $M$, then $\beta^*(A) = 2\beta(A)$.
\end{theorem}
\begin{proof}
    Every $\ep$-diameter ball cover is a $\ep$-diameter set cover, so therefore
    $\alpha(A) \leq \beta^*(A)$.  Hence by Theorem~\ref{T:a=2bmt},
    $2\beta(A) \leq \beta^*(A)$.

    Conversely since every $\ep$-radius ball cover is a $2\ep$-diameter ball cover,
    $\beta^*(A) \leq 2\beta(A)$.  Therefore $\beta^*(A) = 2\beta(A)$.
\end{proof}

\begin{definition}\label{D:setballcontract}
    Let $T:X \to Y$ be a continuous map between metric spaces.

    For $k \geq 0$, $T$ is \emph{$k$-set-contractive} iff for every bounded subset
    $A \subset X$, that $\alpha(T(A)) \leq k\,\alpha(A)$.  Analogously, $T$ is
    \emph{$k$-ball-contractive} iff for all bounded subsets $A \subset X$,
    $\beta(T(A)) \leq k\,\beta(A)$.

    $T$ is \emph{set-condensing (ball-condensing)} iff
    $T$ is $k$-set-contractive ($k$-ball contractive) with $k <1$.
\end{definition}

\begin{lemma}\label{L:setballcont in ms}
    Let $T:X \to Y$ be a continuous map between metric spaces.
    If $T$ is $k$-set-contractive, then $T$ is $2k$-ball-contractive.
    If $T$ is $k$-ball-contractive, then $T$ is $2k$-set-contractive.
\end{lemma}
\begin{proof}
    Suppose $T$ is $k$-set-contractive and let $A \subset X$ be bounded.
    Then by Lemma~\ref{L:ba2b}, we know the following:
    \[
    \beta(T(A)) \leq \alpha(T(A)) \leq k\,\alpha(A) \leq 2k\, \beta(A).
    \]
    Therefore, $T$ is $2k$-ball-contractive.

    Analogously, suppose $T$ is $k$-ball-contractive and let $A \subset X$ be
    bounded.    Then by Lemma~\ref{L:ba2b}, we know the following:
    \[
    \alpha(T(A)) \leq 2\,\beta(T(A)) \leq 2k\,\beta(A) \leq 2k\,\alpha(A).
    \]
    Therefore $T$ is $2k$-set-contractive.
\end{proof}

\begin{theorem}\label{T:setballcont in mt}
    Let $M, N$ be metric trees and $D \subset M$.  Then $T$ is $k$-set-contractive
    iff $T$ is $k$-ball-contractive.
\end{theorem}
\begin{proof}
    Let $A \subset D$ be bounded.

    Suppose that $T$ is $k$-set-contractive, then by Theorem~\ref{T:a=2bmt},
    \[
    \beta(T(A)) = \frac{1}{2}\,\alpha(T(A)) \leq \frac{k}{2}\,\alpha(A) = k\,\beta(A).
    \]
    Therefore $T$ is $k$-ball-contractive.

    Suppose that $T$ is $k$-ball-contractive, then by Theorem~\ref{T:a=2bmt},
    \[
    \alpha(T(A)) = 2\,\beta(T(A)) \leq 2k\, \beta(A) = k\,\alpha(A).
    \]
    Therefore $T$ is $k$-set-contractive.
\end{proof}

We will now examine the Lifschitz characteristic of a metric space
and its connection to measures of noncompactness, (see \cite{Kras}).
\begin{definition}\label{D:lif}
    For a metric space $X$, $b \in \R_{>0}$ is \emph{Lifschitz for $X$}
    iff there exists $a>1$ such that for all $x,y \in X, r>0$,
    if $xy>r$, then there exists $z \in X$ such that
    $B_c(x;ar) \cap B_c(y;br) \subset B_c(z;r)$.
\end{definition}

\begin{definition}\label{D:lifc}
    For a metric space $X$, define the
    \emph{Lifschitz characteristic of $X$} as
    \[
    \kappa(X) := \sup \{ b>0 \mid b \mbox{ is Lifschitz for } X\}.
    \]
\end{definition}

\begin{remark}\label{R:k >= 1}
    Observe that for any metric space $X$, that $\kappa(X) \geq 1$.
    For if $b \leq 1$, then by choosing $z=y$ we can see for any
    $a>1$ that, $B_c(x;ar) \cap B_c(y;br) \subset B_c(z;r)$.
\end{remark}

\begin{theorem}\label{T:kM=2}
    Let $M$ be a metric tree, then $\kappa(M) = 2$.
\end{theorem}
\begin{proof}[Proof of $\kappa \geq 2$]
    We will show that every $0<b<2$ is Lifschitz for $M$.
    Let $b = 2-2\ep$, where $0< \ep < 1$ and select $a>1$ to be
    $a = 1 + \ep$.      For every $r>0$ and $x,y \in M$ with $xy>r$,
    let $z \in [x,y]$ such that $xz = \ep r$.

    Suppose that $w \in B_c(x;ar) \cap B_c(y;br)$.
    By Theorem~\ref{T:msip}, there exists $u \in [x,y]$ such that
    $[w,x] \cap [w,y] = [w,u]$.  Since $z \in [x,y]$ too, we know that
    $z \in [x,u]$ or $z \in [u,y]$.

    If $z \in [x,u]$, then by $u \in [x,w]$ and transitivity, we know
    $z \in [x,w]$.  Hence,
$$
\begin{array}{rcll}
 zw &=& xw - xz \leq ar - xz  &\mbox{ by } z \in [x,w] \mbox{ and } w \in B_c(x;ar), \\
      &=& (r + \ep r) - \ep r &\mbox{ by } a = 1+\ep \mbox{ and } xz =\ep r, \\
      &=& r. &
\end{array}
$$

    If $z \in [u,y]$, then by $u \in [w,y]$ and transitivity, we know
    $z \in [w,y]$.  Hence,
$$
    \begin{array}{rcll}
    zw &=& wy - zy &\mbox{ by } z \in [w,y]\\
    	&=& wy - (xy - xz) &\mbox{ by } z \in [x,y],\\
        &<& br - r + \ep r &\mbox{ by } w \in B_c(y;br),\, r < xy, \mbox{ and } xz = \ep r,\\
        &=& (2r - 2\ep r) - r + \ep r &\mbox{ by } b=2-2\ep,\\
        &=& r - \ep r < r. &
   \end{array}
$$
    We know that $z \in [x,u]$ or $z \in [u,y]$, and in either case $zw \leq r$, so
    $w \in B_c(z; r)$.  Therefore $B_c(x;ar) \cap B_c(y;br) \subset B_c(z;r)$.  
    Hence $b$ is
    Lifschitz for $M$ for all $b<2$ and therefore $\kappa(M) \geq 2$.
\end{proof}
\begin{proof}[Proof of $\kappa(M) \leq 2$]
    Suppose for some metric tree $M$, that $b$ is Lifschitz for $M$.
    Hence there exists $a>1$ such that for all $r>0$ and $x,y \in M$
    with $xy>r$, there exists $z \in M$ such that
    $B_c(x;ar) \cap B_c(y;2r) \subset B_c(z;r)$.

    Let $w,v \in M$ with $wv = 4r$ for some $r>0$ and let $a>1$.
    Let $y$ be the midpoint of $[w,v]$ and let $x \in [y,v]$ with
    $r < yx < ar$.  Hence $[w,v] \subset B_c(y;2r)$ and
    $[u,v] \subset B_c(x;ar)$ where $u \in [w,v]$ with $ux = ar$.
    Therefore, $[u,v] \subset B_c(x;ar) \cap B_c(y;2r)$.

    However, $diam([u,v]) = uy + yv > yv = 2r$.  So for all $z \in M$,
    $[u,v] \nsubseteq B_c(z;r)$ since $diam(B_c(z;r)) \leq 2r$.
    Therefore for all 
    $z \in M$, $B_c(x;ar) \cap B_c(y;2r) \nsubseteq B_c(z;r)$, which
    contradicts that $b=2$ is Lifschitz for $M$.  Therefore, for all metric
    trees $M$, $\kappa(M) \leq 2$.
\end{proof}
\begin{lemma}[Webb and Zhao, \cite{Webb}]\label{L:webb}
    Let $M$ be a metric space and $A \subset M$, with $diam(A)=d$.
    For all $0<b<\kappa(M)$, there exists $z \in M$ such that
    $A \subset B_c(z,\frac{d}{b})$.
\end{lemma}

\begin{theorem}[Webb and Zhao, \cite{Webb}]\label{T:webb}
    Let $A$ be a bounded subset of metric space $M$,
    then $\kappa(M)\beta(A) \leq \alpha(A)$.
\end{theorem}
\begin{proof}
    By the definition of $\alpha(A)$, for all $\delta$ with $0<\delta<\kappa(M)$,
    there exists $E_j \subset M$ such that $diam(E_j)\leq \alpha(A) + \delta$
    and $A \subset \bigcup_{j=1}^{n} E_j$.

    By Lemma~\ref{L:webb}, for $b= \kappa(M)-\delta$ and each $E_j$, there
    exists $z_j \in M$ such that
    $E_j \subset B_c(z_j,\frac{\alpha(A)+\delta}{\kappa(M) - \delta})$. Hence,
    \[
    A \subset \bigcup_{j=1}^{n} E_j \subset \bigcup_{j=1}^{n}
    B_c\left(z_j,\frac{\alpha(A)+\delta}{\kappa(M) - \delta}\right).
    \]
    So for all $\delta$ such that $0<\delta<\kappa(M)$, we have that
    $\beta(A) \leq \frac{\alpha(A)+\delta}{\kappa(M) - \delta}$.
    Therefore, $\kappa(M)\beta(A) \leq \alpha(A)$.
\end{proof}

\begin{remark}\label{R:a=2b}
    Theorems~\ref{T:kM=2} and \ref{T:webb}, show that
    $2\beta(A)\leq\alpha(A)$ for any bounded subset $A$ in a metric tree $M$.
    Therefore we have an alternate proof that for any bounded subset $A$
    in a metric tree $M$, $\alpha(A)=2\beta(A)$. Moreover, by
    Theorems 3.4 and 3.5  we know that $\alpha$ is invariant under isometric
    embedding $f:X\rightarrow Y$ and $\beta$ is not. However if $X$
    and $Y$ are metric trees, $\beta$ will be invariant under $f$ as
    well, since $\alpha(A)=2\beta(A)$.
\end{remark}

\begin{IsbfsBIBLIOGRAPHY}

    \bibitem{Akso}
        A.G. Aksoy and M. A. Khamsi, \textit{A Selection Theorem in Metric Trees}, 
        Proc.\ Amer.\ Math.\ Soc.\ {\bf 134}
        (2006), 2957--2966.
    \bibitem{Bana}
        J. Banas and K. Goebel, \textit{Measures of Noncompactness in Banach Spaces},
        Lecture Notes in Pure and Appl. Math, vol. 60, Dekker, New York, 1980.
    \bibitem{Bart}
        I. Bartolini, P. Ciaccia, and M. Patella, \textit{String Matching with Metric Trees
        Using Approximate Distance}, SPIR, Lecture Notes in Computer Science,
        Springer-Verlag, vol. 2476, 2002, 271--283.
    \bibitem{Best}
        M. Bestvina, \textit{$\mathbb{R}$-trees in Topology, Geometry, and Group Theory},
        Handbook of geometric topology, North-Holland, Amsterdam, 2002, 55--91.
    \bibitem{Blum}
        L. M. Blumenthal, \textit{Theory and Applications of Distance Geometry},
        Oxford University Press, London, 1953.
    \bibitem{Brid}
        M. Bridson and A. Haefliger, \textit{Metric Spaces of Nonpositive Curvature},
        Grundlehren der Mathematischen Wissenschaften, vol. 319,
        Springer-Verlag, Berlin, 1999.
    \bibitem{Darb}
        G. Darbo, \textit{Punti Uniti in Trasformazioni a Codominio Non Compatto},
        Rend. Sem. Mat. Univ. Padova {\bf 24} (1955), 84--92.
     \bibitem{Kirk1}
        R. Espinola and W. A. Kirk, \textit{Fixed Point Theorems in $\mathbb{R}$-trees
        with Applications to Graph Theory}, Topology Appl. {\bf 153} (2006),
        1046--1055.
     \bibitem{Kirk2}
        W. A. Kirk, \textit{Fixed Point Theorems in $\rm CAT(0)$ Spaces and
        $\mathbb{R}$-Trees}, Fixed Point Theory Appl. (2004), 309--316.
    \bibitem{Kras}
        M. A. Krasnosel'skii and P. P. Zabreiko, \textit{Geometric Methods of Nonlinear
        Analysis}, Springer, Berlin, 1984.
    \bibitem{Kura}
        K. Kuratowski, \textit{Sur les Espaces Complets}, Fund. Math.
        {\bf 15} (1930), 301--309.
    \bibitem{Semp}
        C. Semple and M. Steel, \textit{Phylogenetics}, Oxford Lecture Series
        in Mathematics and its Applications, vol. 24, 2003.
    \bibitem{Tits}
        J. Tits, \textit{A Theorem of Lie-Kolchin for Trees},
        Contributions to Algebra: A Collection of Papers Dedicated to Ellis Kolchin,
        Academic Press, New York, 1977.
    \bibitem{Webb}
        J.R.L. Webb and W. Zhao, \textit{On Connections Between Set and Ball
        Measures of Noncompactness}, Bull. London Math. Soc. {\bf 22} (1990),
        471--477.
\end{IsbfsBIBLIOGRAPHY}

\end{IsbfsDOCUMENT}